\documentclass{article}
\usepackage{amssymb}
\usepackage{amsfonts}
\usepackage{amsmath}
\usepackage{harvard}

\newtheorem{theorem}{Theorem}

\newtheorem{corollary}[theorem]{Corollary}

\newtheorem{definition}[theorem]{Definition}
\newtheorem{example}[theorem]{Example}

\newtheorem{proposition}[theorem]{Proposition}
\newtheorem{remark}[theorem]{Remark}

\newenvironment{proof}[1][Proof]{\noindent\textbf{#1.} }{\ \rule{0.5em}{0.5em}}
\numberwithin{theorem}{section}
\numberwithin{equation}{section}
\input{tcilatex}

\begin{document}

\title{Deflection d-Tensor Identities in the Relativistic Time Dependent
Lagrange Geometry}
\author{Mircea Neagu and Emil Stoica}
\date{}
\maketitle

\begin{abstract}
The aim of this paper is to study the local components of the relativistic
time dependent d-linear connections, d-torsions, d-curvatures and deflection
d-tensors with respect to an adapted basis on the 1-jet space $J^{1}(\mathbb{%
R},M)$. The Ricci identities, together with their corresponding identities
of deflection d-tensors, are also given.
\end{abstract}



\section{Some physical and geometrical aspects}

\hspace{5mm}According to Olver's opinion [10], we agree that the 1-jet fibre
bundle is a basic object in the study of classical and quantum field
theories. For such a reason, a lot of authors (Asanov [2], Saunders [11],
Vondra [12] and many others) studied the differential geometry of the 1-jet
spaces. Continuing the geometrical studies of Asanov [2] and using as a
pattern the Lagrangian geometrical ideas developed by Miron, Anastasiei and
Buc\u{a}taru in the monographs [5] and [3], the first author of this paper
has recently developed the \textit{Riemann-Lagrange geometry of 1-jet spaces}
[7]. This theory is very suitable for the geometrical study of the \textit{%
relativistic non-autonomous (rheonomic) Lagrangians}, that is of the
Lagrangians depending on an usual \textit{relativistic time} [6], [8] or
depending on a \textit{relativistic multi-time} [7], [9].

It is important to note that a \textit{classical non-autonomous (rheonomic)
Lagrangian geometry} (i. e. a geometrization of the Lagrangians depending on
an \textit{absolute time}) was sketched by Miron and Anastasiei at the end
of the book [5] and developed in the same way by Anastasiei and Kawaguchi
[1] or Frigioiu [4].

In what follows we try to expose the main geometrical and physical aspects
which differentiate the both geometrical theories: the \textit{jet
relativistic non-autonomous Lagrangian geometry} [8] and the \textit{%
classical non-a\-u\-to\-no\-mous Lagrangian geometry} [5].

In this direction, we point out that the \textit{relativistic non-autonomous
Lagrangian geometry} [8] has as natural house the 1-jet space $J^{1}(\mathbb{%
R},M)$, where $\mathbb{R}$ is the manifold of real numbers having the
coordinate $t$. This represents for us the usual \textit{relativistic time}.
We recall that the 1-jet space $J^{1}(\mathbb{R},M)$ is regarded as a vector
bundle over the product manifold $\mathbb{R}\times M$, having the fibre type 
$\mathbb{R}^{n}$, where $n$ is the dimension of the \textit{spatial}
manifold $M$. In mechanical terms, if the manifold $M$ has the spatial local
coordinates $(x^{i})_{i=\overline{1,n}}$, then the 1-jet vector bundle%
\begin{equation*}
J^{1}(\mathbb{R},M)\rightarrow \mathbb{R}\times M
\end{equation*}%
can be regarded as a \textit{bundle of configurations} having the local
coordinates $(t,x^{i},y_{1}^{i})$; these transform by the rules [8]%
\begin{equation}
\left\{ 
\begin{array}{l}
\widetilde{t}=\widetilde{t}(t)\medskip \\ 
\widetilde{x}^{i}=\widetilde{x}^{i}(x^{j})\medskip \\ 
\widetilde{y}_{1}^{i}=\dfrac{\partial \widetilde{x}^{i}}{\partial x^{j}}%
\dfrac{dt}{d\widetilde{t}}\cdot y_{1}^{j}.%
\end{array}%
\right.  \label{rgg}
\end{equation}

\begin{remark}
The form of the jet transformation group (\ref{rgg}) stands out by the 
\textbf{relativistic }\textit{character} of the \textbf{time} $t$.
\end{remark}

Comparatively, in the \textit{classical non-a\-u\-to\-no\-mous Lagrangian
geometry} [5] the \textit{bundle of configurations }is the vector bundle%
\begin{equation*}
\mathbb{R}\times TM\rightarrow M,
\end{equation*}%
whose local coordinates $(t,x^{i},y^{i})$ transform by the rules 
\begin{equation}
\left\{ 
\begin{array}{l}
\widetilde{t}=t\medskip  \\ 
\widetilde{x}^{i}=\widetilde{x}^{i}(x^{j})\medskip  \\ 
\widetilde{y}^{i}=\dfrac{\partial \widetilde{x}^{i}}{\partial x^{j}}\cdot
y^{j},%
\end{array}%
\right.   \label{agg}
\end{equation}%
where $TM$ is the tangent bundle of the spatial manifold $M$.

\begin{remark}
The form of the transformation group (\ref{agg}) stands out by the \textbf{%
absolute }\textit{character} of the \textbf{time} $t$.
\end{remark}

It is important to note that jet transformation group (\ref{rgg}) from the 
\textit{relativistic non-autonomous Lagrangian geometry} is more general and
more natural than the transformation group (\ref{agg}) used in the \textit{%
classical non-autonomous Lagrangian geometry}. This is because the last one
ignores the temporal reparametrizations, emphasizing in this way the
absolute character of the usual time coordinate $t$. Or, physically
speaking, the relativity of time is an well-known fact.

From a geometrical point of view, we point out that the entire \textit{%
classical rheonomic Lagrangian geometry} of Miron and Anastasiei [5] relies
on the study of the \textit{energy action functional}%
\begin{equation*}
\mathbb{E}_{1}(c)=\int_{a}^{b}L(t,x^{i},y^{i})dt,
\end{equation*}%
where $L:\mathbb{R}\times TM\rightarrow \mathbb{R}$ is a Lagrangian function
and $y^{i}=dx^{i}/dt,$ whose Euler-Lagrange equations produce a semispray $%
G^{i}(t,x^{k},y^{k})$ and a corresponding nonlinear connection $%
N_{j}^{i}=\partial G^{i}/\partial y^{j}.$ Therefore, the authors construct
the adapted bases of vector and covector fields, together with the adapted
components of the $N$-linear connections and their corresponding d-torsions
and d-curvatures. But, because $L(t,x^{i},y^{i})$ is a real function, we
deduce that the previous geometrical theory has the following impediment: -%
\textit{the energy action functional depends on the reparametrizations }$%
t\longleftrightarrow \widetilde{t}$\textit{\ of the same curve }$c.$ Thus,
in order to avoid this inconvenience, the Finsler case imposes the
1-positive homogeneity condition%
\begin{equation*}
L(t,x^{i},\lambda y^{i})=\lambda L(t,x^{i},y^{i}),\text{ }\forall \text{ }%
\lambda >0.
\end{equation*}

Alternatively, the \textit{relativistic rheonomic Lagrangian geometry} from
[8] uses the \textit{relativistic energy action functional}%
\begin{equation*}
\mathbb{E}_{2}(c)=\int_{a}^{b}L(t,x^{i},y_{1}^{i})\sqrt{h_{11}(t)}dt,
\end{equation*}%
where $L:J^{1}(\mathbb{R},M)\rightarrow \mathbb{R}$ is a jet Lagrangian
function and $h_{11}(t)$ is a Riemannian metric on the relativistic time
manifold $\mathbb{R}$. This functional is now independent by the
reparametrizations $t\longleftrightarrow \widetilde{t}$ of the same curve $c$%
. The Euler-Lagrange equations of the Lagrangian $\mathcal{L}%
=L(t,x^{i},y_{1}^{i})\sqrt{h_{11}(t)}$ produce a relativistic time dependent
semispray [8]%
\begin{equation*}
\mathcal{S=}\left( H_{(1)1}^{(i)},\text{ }G_{(1)1}^{(i)}\right) ,
\end{equation*}%
which gives the jet nonlinear connection [6]%
\begin{equation*}
\Gamma _{\mathcal{S}}=\left( M_{(1)1}^{(j)}=2H_{(1)1}^{(j)},\text{ }%
N_{(1)k}^{(j)}=\frac{\partial G_{(1)1}^{(j)}}{\partial y_{1}^{k}}\right) .
\end{equation*}%
With these geometrical tools we can construct in the relativistic rheonomic
Lagrangian geometry the distinguished (d-) linear connections, together with
their d-torsions and d-curvatures, which naturally generalize the similar
geometrical objects from the classical rheonomic Lagrangian geometry [5]. In
this respect, the authors of this paper believe that the relativistic
geometrical approach proposed in this paper has more geometrical and
physical meanings than the theory proposed by Miron and Anastasiei in [5].

In conclusion, in order to remark the main similitudes and differences
between these geometrical theories, we invite the reader to compare both the 
\textit{classical }and \textit{relativistic non-autonomous Lagrangian
geometries} exposed in the works [5] and [8].

As a final remark, we point out that for a lot of mathematicians (such as
Crampin, de Leon, Krupkova, Sarlet, Saunders and others) the non-autonomous
Lagrangian geometry is constructed on the first jet bundle $J^{1}\pi $ of a
fibered manifold $\pi :M^{n+1}\longrightarrow \mathbb{R}.$ In their papers,
if $(t,x^{i})$ are the local coordinates on the $n+1$-dimensional manifold $M
$ such that $t$ is a global coordinate for the fibers of the submersion $\pi 
$ and $x^{i}$ are transverse coordinates of the induced foliation, then a
change of coordinates on $M$ is given by%
\begin{equation}
\left\{ 
\begin{array}{ll}
\widetilde{t}=\widetilde{t}(t),\medskip  & \dfrac{d\widetilde{t}}{dt}\neq 0
\\ 
\widetilde{x}^{i}=\widetilde{x}^{i}(x^{j},t), & \text{rank}\left( \dfrac{%
\partial \widetilde{x}^{i}}{\partial x^{j}}\right) =n.%
\end{array}%
\right.   \label{Krupkova}
\end{equation}

Altough the 1-jet extension of the transformation rules (\ref{Krupkova}) is
more general than the transformation group (\ref{rgg}), the authors ot this
paper consider that the transformation group (\ref{rgg}) is more appropriate
for their final purpose, the development of a \textit{relativistic rheonomic
Lagrangian field theory}. For example, in the paper [8], starting with a
non-degenerate Lagrangian function $L:J^{1}(\mathbb{R},M)\rightarrow \mathbb{%
R}$ and an \textit{a priori} given Riemannian metric $h_{11}(t)$ on the
relativistic temporal manifold $\mathbb{R}$, one introduces a \textit{%
relativistic time dependent electromagnetic field}%
\begin{equation*}
F=F_{(i)j}^{(1)}\delta y_{1}^{i}\wedge dx^{j},
\end{equation*}%
where%
\begin{equation*}
F_{(i)j}^{(1)}=\frac{1}{2}\left[ D_{(i)j}^{(1)}-D_{(j)i}^{(1)}\right] \text{
and }\delta y_{1}^{i}=dy_{1}^{i}+M_{(1)1}^{(i)}dt+N_{(1)j}^{(i)}dx^{j},
\end{equation*}%
the \textit{metrical deflection d-tensors} $D_{(i)j}^{(1)}$ beeing produced
only by the jet Lagrangian $\mathcal{L}=L\sqrt{h_{11}(t)}.$ In such a
perspective, the relativistic time dependent electromagnetic field $F$ has
an \textit{intrinsic geometrical character}. Moreover, the electromagnetic
components $F_{(i)j}^{(1)}$ are governed by some natural \textit{generalized
Maxwell equations}. These equations are exposed in [8] and naturally
generalize the already classical Maxwell equations from Miron and
Anastasiei's theory [5].

\section{The adapted components of the jet $\Gamma $-linear connections}

\label{coli}

\hspace{5mm}Let us suppose that on the 1-jet space $E=J^{1}(\mathbb{R},M)$
is fixed a nonlinear connection $\Gamma $ given by the \textit{temporal
components} $M_{(1)1}^{(i)}$ and the \textit{spatial components} $%
N_{(1)j}^{(i)}$. We recall that the transformation rules of the local
components of the nonlinear connection $\Gamma =\left(
M_{(1)1}^{(i)},N_{(1)j}^{(i)}\right) $ are expressed by [6]%
\begin{equation*}
\widetilde{M}_{(1)1}^{(k)}=M_{(1)1}^{(j)}\left( \dfrac{dt}{d\widetilde{t}}%
\right) ^{2}\dfrac{\partial \widetilde{x}^{k}}{\partial x^{j}}-\dfrac{dt}{d%
\widetilde{t}}\dfrac{\partial \widetilde{y}_{1}^{k}}{\partial t},\quad%
\widetilde{N}_{(1)l}^{(k)}=N_{(1)i}^{(j)}\dfrac{dt}{d\widetilde{t}}\dfrac{%
\partial x^{i}}{\partial \widetilde{x}^{l}}\dfrac{\partial \widetilde{x}^{k}%
}{\partial x^{j}}-\dfrac{\partial x^{i}}{\partial \widetilde{x}^{l}}\dfrac{%
\partial \widetilde{y}_{1}^{k}}{\partial x^{i}}.
\end{equation*}

\begin{example}
\label{cannlc} Let $(\mathbb{R},h_{11}(t))$ and $(M,\varphi _{ij}(x))$ be
semi-Rie\-ma\-nni\-an manifolds. Let us consider the Christoffel symbols%
\begin{equation*}
H_{11}^{1}=\frac{h^{11}}{2}\frac{dh_{11}}{dt}
\end{equation*}%
and%
\begin{equation*}
\gamma _{jk}^{i}=\frac{\varphi ^{im}}{2}\left( \frac{\partial \varphi _{jm}}{%
\partial x^{k}}+\frac{\partial \varphi _{km}}{\partial x^{j}}-\frac{\partial
\varphi _{jk}}{\partial x^{m}}\right) ,
\end{equation*}%
where $h^{11}=1/h_{11}$. Then, using the transformation rules%
\begin{equation}
\widetilde{H}_{11}^{1}=H_{11}^{1}\frac{dt}{d\widetilde{t}}+\frac{d\widetilde{%
t}}{dt}\frac{d^{2}t}{d\widetilde{t}^{2}}  \label{tr-rules-H}
\end{equation}%
and%
\begin{equation}
\widetilde{\gamma }_{qr}^{p}=\gamma _{jk}^{i}\frac{\partial \widetilde{x}^{p}%
}{\partial x^{i}}\frac{\partial x^{j}}{\partial \widetilde{x}^{q}}\frac{%
\partial x^{k}}{\partial \widetilde{x}^{r}}+\frac{\partial \widetilde{x}^{p}%
}{\partial x^{l}}\frac{\partial ^{2}x^{l}}{\partial \widetilde{x}%
^{q}\partial \widetilde{x}^{r}},  \label{tr-rules-gamma}
\end{equation}%
we deduce that the set of local functions%
\begin{equation*}
\mathring{\Gamma}=\left( \mathring{M}_{(1)1}^{(j)},\mathring{N}%
_{(1)i}^{(j)}\right) ,
\end{equation*}%
where%
\begin{equation}
\mathring{M}_{(1)1}^{(j)}=-H_{11}^{1}y_{1}^{j}\text{ \ \ and \ \ }\mathring{N%
}_{(1)i}^{(j)}=\gamma _{im}^{j}y_{1}^{m},  \label{nlc-assoc-metrics}
\end{equation}%
represents a nonlinear connection on the 1-jet space $E=J^{1}(\mathbb{R},M)$%
. This jet nonlinear connection is called the \textbf{canonical nonlinear
connection attached to the pair of metrics} $(h_{11}(t),\varphi _{ij}(x))$.
\end{example}

Let 
\begin{equation*}
{\left\{ {\frac{\delta }{\delta t}},{\frac{\delta }{\delta x^{i}}},{\frac{%
\partial }{\partial y_{1}^{i}}}\right\} \subset \mathcal{X}(E)}
\end{equation*}%
and%
\begin{equation*}
\{dt,dx^{i},\delta y_{1}^{i}\}\subset \mathcal{X}^{\ast }(E)
\end{equation*}%
be the dual bases adapted to the nonlinear connection $\Gamma =\left(
M_{(1)1}^{(i)},N_{(1)j}^{(i)}\right) ,$ where%
\begin{equation*}
\begin{array}{l}
{{\dfrac{\delta }{\delta t}}={\dfrac{\partial }{\partial t}}-M_{(1)1}^{(j)}{%
\dfrac{\partial }{\partial y_{1}^{j}}}},\medskip \\ 
{{\dfrac{\delta }{\delta x^{i}}}={\dfrac{\partial }{\partial x^{i}}}%
-N_{(1)i}^{(j)}{\dfrac{\partial }{\partial y_{1}^{j}}}},\medskip  \\ 
\delta y_{1}^{i}=dy_{1}^{i}+M_{(1)1}^{(i)}dt+N_{(1)j}^{(i)}dx^{j}.%
\end{array}%
\end{equation*}

\begin{remark}
The components of the above dual adapted bases transform under a change of
coordinates (\ref{rgg}) as classical tensors.
\end{remark}

In order to develop a theory of the $\Gamma $-linear connections on the
1-jet space $E=J^{1}(\mathbb{R},M)$, we need the following simple result:

\begin{proposition}
\textbf{(a)} The Lie algebra $\mathcal{X}(E)$ of the vector fields on $E$
decomposes in the direct sum%
\begin{equation*}
\mathcal{X}(E)=\mathcal{X}(\mathcal{H}_{\mathbb{R}})\oplus \mathcal{X}(%
\mathcal{H}_{M})\oplus \mathcal{X}(\mathcal{V}),
\end{equation*}%
where%
\begin{equation*}
\mathcal{X}(\mathcal{H}_{\mathbb{R}})=Span\left\{ {\frac{\delta }{\delta t}}%
\right\} ,\quad \mathcal{X}(\mathcal{H}_{M})=Span\left\{ {\frac{\delta }{%
\delta x^{i}}}\right\} ,\quad \mathcal{X}(\mathcal{V})=Span\left\{ {\frac{%
\partial }{\partial y_{1}^{i}}}\right\} .
\end{equation*}

\textbf{(b)} The Lie algebra $\mathcal{X}^{\ast }(E)$ of the covector fields
on $E$ decomposes in the direct sum%
\begin{equation*}
\mathcal{X}^{\ast }(E)=\mathcal{X}^{\ast }(\mathcal{H}_{\mathbb{R}})\oplus 
\mathcal{X}^{\ast }(\mathcal{H}_{M})\oplus \mathcal{X}^{\ast }(\mathcal{V}),
\end{equation*}%
where%
\begin{equation*}
\mathcal{X}^{\ast }(\mathcal{H}_{\mathbb{R}})=Span\{dt\},\quad \mathcal{X}%
^{\ast }(\mathcal{H}_{M})=Span\{dx^{i}\},\quad \mathcal{X}^{\ast }(\mathcal{V%
})=Span\{\delta y_{1}^{i}\}.
\end{equation*}
\end{proposition}

Denoting by $h_{\mathbb{R}}$, $h_{M}$, respectively $v$, the $\mathbb{R}$%
\textit{-horizontal}, $M$\textit{-horizontal}, respectively \textit{vertical
canonical projections} associated to the above decompositions, we get

\begin{corollary}
\textbf{(a)} Any vector field on $E$ can be uniquely written in the form: 
\begin{equation*}
X=h_{\mathbb{R}}X+h_{M}X+vX,\quad \forall \text{ }X\in \mathcal{X}(E).
\end{equation*}

\textbf{(b)} Any 1-form on $E$ can be uniquely written in the form: 
\begin{equation*}
\omega =h_{\mathbb{R}}\omega +h_{M}\omega +v\omega ,\quad \forall \text{ }%
\omega \in \mathcal{X}^{\ast }(E).
\end{equation*}
\end{corollary}

\begin{definition}
A linear connection $\nabla :\mathcal{X}(E)\times \mathcal{X}(E)\rightarrow 
\mathcal{X}(E),$ which verifies the Ehresman-Koszul axioms 
\begin{equation*}
\nabla h_{\mathbb{R}}=0,\quad \nabla h_{M}=0,\quad \nabla v=0,
\end{equation*}%
is called a $\Gamma $\textbf{-linear connection}\textit{\ on the 1-jet
vector bundle} $E=J^{1}(\mathbb{R},M)$.
\end{definition}

Using the adapted basis of vector fields on $E=J^{1}(\mathbb{R},M)$ and the
definition of a $\Gamma $-linear connection, we prove without difficulties

\begin{proposition}
A $\Gamma $-linear connection $\nabla $ on $E=J^{1}(\mathbb{R},M)$ is
determined by \textbf{nine} local adapted components%
\begin{equation*}
\nabla \Gamma =\left( \bar{G}_{11}^{1},\;G_{i1}^{k},\;G_{(1)(j)1}^{(i)(1)},\;%
\bar{L}_{1j}^{1},\;L_{ij}^{k},\;L_{(1)(j)k}^{(i)(1)},\;\bar{C}%
_{1(k)}^{1(1)},\;C_{i(k)}^{j(1)},\;C_{(1)(j)(k)}^{(i)(1)(1)}\right) ,
\end{equation*}%
which are uniquely defined by the relations:%
\begin{equation*}
\begin{array}{l}
(h_{\mathbb{R}})\hspace{4mm}{\nabla _{{\dfrac{\delta }{\delta t}}}{\dfrac{%
\delta }{\delta t}}=\bar{G}_{11}^{1}{\dfrac{\delta }{\delta t}},\;\;\nabla _{%
{\dfrac{\delta }{\delta t}}}{\dfrac{\delta }{\delta x^{i}}}=G_{i1}^{k}{%
\dfrac{\delta }{\delta x^{k}}},\;\;\nabla _{{\dfrac{\delta }{\delta t}}}{%
\dfrac{\partial }{\partial y_{1}^{i}}}=G_{(1)(i)1}^{(k)(1)}{\dfrac{\partial 
}{\partial y_{1}^{k}}}},\medskip  \\ 
(h_{M})\;\;{\nabla _{{\dfrac{\delta }{\delta x^{j}}}}{\dfrac{\delta }{\delta
t}}=\bar{L}_{1j}^{1}{\dfrac{\delta }{\delta t}},\;\;\nabla _{{\dfrac{\delta 
}{\delta x^{j}}}}{\dfrac{\delta }{\delta x^{i}}}=L_{ij}^{k}{\dfrac{\delta }{%
\delta x^{k}}},\;\;\nabla _{{\dfrac{\delta }{\delta x^{j}}}}{\dfrac{\partial 
}{\partial y_{1}^{i}}}=L_{(1)(i)j}^{(k)(1)}{\dfrac{\partial }{\partial
y_{1}^{k}}}},\medskip  \\ 
(v)\hspace{3mm}{\nabla _{{\dfrac{\partial }{\partial y_{1}^{j}}}}{\dfrac{%
\delta }{\delta t}}=\bar{C}_{1(j)}^{1(1)}{\dfrac{\delta }{\delta t}}%
,\;\nabla _{{\dfrac{\partial }{\partial y_{1}^{j}}}}{\dfrac{\delta }{\delta
x^{i}}}=C_{i(j)}^{k(1)}{\dfrac{\delta }{\delta x^{k}}},\;\nabla _{{\dfrac{%
\partial }{\partial y_{1}^{j}}}}{\dfrac{\partial }{\partial y_{1}^{i}}}%
=C_{(1)(i)(j)}^{(k)(1)(1)}{\dfrac{\partial }{\partial y_{1}^{k}}}}.%
\end{array}%
\end{equation*}
\end{proposition}

Taking into account the tensorial transformation laws of the adapted basis
of vector fields on $E=J^{1}(\mathbb{R},M)$, by laborious local
computations, we deduce

\begin{theorem}
\textbf{(a)} Under a change of coordinates (\ref{rgg}) on the 1-jet vector
bundle $E=J^{1}(\mathbb{R},M)$, the adapted coefficients of the $\Gamma $%
-linear connection $\nabla $ modify by the rules\medskip \newline
$(h_{\mathbb{R}})\hspace{6mm}\left\{ 
\begin{array}{l}
\medskip {\bar{G}_{11}^{1}=\tilde{\bar{G}}_{11}^{1}{\dfrac{d\tilde{t}}{dt}}+{%
\dfrac{dt}{d\tilde{t}}\dfrac{d^{2}\tilde{t}}{dt^{2}}}} \\ 
\medskip {G_{i1}^{k}=\tilde{G}_{j1}^{r}{\dfrac{\partial x^{k}}{\partial 
\tilde{x}^{r}}}{\dfrac{\partial \tilde{x}^{j}}{\partial x^{i}}}{\dfrac{d%
\tilde{t}}{dt}}} \\ 
{G_{(1)(i)1}^{(k)(1)}=\tilde{G}_{(1)(j)1}^{(p)(1)}{\dfrac{\partial x^{k}}{%
\partial \tilde{x}^{p}}\dfrac{\partial \tilde{x}^{j}}{\partial x^{i}}\dfrac{d%
\tilde{t}}{dt}}+\delta _{i}^{k}}\left( {{\dfrac{d\tilde{t}}{dt}}}\right) ^{2}%
{{\dfrac{d^{2}t}{d\tilde{t}^{2}}}},%
\end{array}%
\right. \medskip $\newline
$(h_{M})\hspace{5mm}\left\{ 
\begin{array}{l}
\medskip {\bar{L}_{1j}^{1}=\tilde{\bar{L}}_{1l}^{1}{\dfrac{\partial \tilde{x}%
^{l}}{\partial x^{j}}}} \\ 
\medskip {L_{ij}^{r}=\tilde{L}_{pq}^{s}{\dfrac{\partial x^{r}}{\partial 
\tilde{x}^{s}}\dfrac{\partial \tilde{x}^{p}}{\partial x^{i}}}{\dfrac{%
\partial \tilde{x}^{q}}{\partial x^{j}}}+{\dfrac{\partial x^{r}}{\partial 
\tilde{x}^{s}}\dfrac{\partial ^{2}\tilde{x}^{s}}{\partial x^{i}\partial x^{j}%
}}} \\ 
{L_{(1)(i)j}^{(k)(1)}=\tilde{L}_{(1)(p)s}^{(r)(1)}{\dfrac{\partial x^{k}}{%
\partial \tilde{x}^{r}}\dfrac{\partial \tilde{x}^{p}}{\partial x^{i}}\dfrac{%
\partial \tilde{x}^{s}}{\partial x^{j}}}+{\dfrac{\partial x^{k}}{\partial 
\tilde{x}^{r}}}{\dfrac{\partial ^{2}\tilde{x}^{r}}{\partial x^{i}\partial
x^{j}}}},%
\end{array}%
\right. \medskip $\newline
$(v)\hspace{8mm}\left\{ 
\begin{array}{l}
\medskip {\bar{C}_{1(i)}^{1(1)}=\tilde{\bar{C}}_{1(j)}^{1(1)}{\dfrac{%
\partial \tilde{x}^{j}}{\partial x^{i}}}{\dfrac{dt}{d\tilde{t}}}} \\ 
\medskip {C_{i(j)}^{k(1)}=\tilde{C}_{p(r)}^{s(1)}{\dfrac{\partial x^{k}}{%
\partial \tilde{x}^{s}}}{\dfrac{\partial \tilde{x}^{p}}{\partial x^{i}}}{%
\dfrac{\partial \tilde{x}^{r}}{\partial x^{j}}}{\dfrac{dt}{d\tilde{t}}}} \\ 
{C_{(1)(i)(j)}^{(k)(1)(1)}=\tilde{C}_{(1)(p)(q)}^{(r)(1)(1)}{\dfrac{\partial
x^{k}}{\partial \tilde{x}^{r}}\dfrac{\partial \tilde{x}^{p}}{\partial x^{i}}%
\dfrac{\partial \tilde{x}^{q}}{\partial x^{j}}}{\dfrac{dt}{d\tilde{t}}}}.%
\end{array}%
\right. \medskip $

\textbf{(b)} Conversely, to give a $\Gamma $-linear connection $\nabla $ on
the 1-jet vector bundle $E=J^{1}(\mathbb{R},M)$ is equivalent to give a set
of \textbf{nine} adapted local components $\nabla \Gamma $, which transform
by the rules described in \textbf{(a)}.
\end{theorem}

\begin{example}
Let $h_{11}(t)$ (respectively $\varphi _{ij}(x)$) be a semi-Riemannian
metric on $\mathbb{R}$ (respectively $M$). We denote by $H_{11}^{1}(t)$
(respectively $\gamma _{ij}^{k}(x)$) the Christoffel symbols of the metric $%
h_{11}(t)$ (respectively $\varphi _{ij}(x)$). Let us consider on the 1-jet
space $E=J^{1}(\mathbb{R},M)$ the canonical nonlinear connection $\mathring{%
\Gamma}$ associated to the pair of metrics $(h_{11}(t),\varphi _{ij}(x))$,
which is defined by the local coefficients (\ref{nlc-assoc-metrics}). In
this context, using the transformation laws (\ref{tr-rules-H}) and (\ref%
{tr-rules-gamma}), we deduce that the set of adapted local coefficients%
\begin{equation*}
B\mathring{\Gamma}=\left( \bar{G}_{11}^{1},\;0,\;G_{(1)(i)1}^{(k)(1)},\;0,%
\;L_{ij}^{k},\;L_{(1)(i)j}^{(k)(1)},\;0,\;0,\;0\right) ,
\end{equation*}%
where%
\begin{equation*}
\bar{G}_{11}^{1}=H_{11}^{1},\;\;G_{(1)(i)1}^{(k)(1)}=-\delta
_{i}^{k}H_{11}^{1},\;\;L_{ij}^{k}=\gamma
_{ij}^{k},\;\;L_{(1)(i)j}^{(k)(1)}=\gamma _{ij}^{k},
\end{equation*}%
defines a $\mathring{\Gamma}$-linear connection on the 1-jet space $E$,
which is called the \textbf{Berwald connection attached to the
semi-Riemannian metrics }$h_{11}(t)$ \textbf{and} $\varphi _{ij}(x)$.
\end{example}

\begin{remark}
In the particular case $(\mathbb{R},h)=(\mathbb{R},\delta )$ our Berwald
linear connection naturally generalizes the canonical $N$-linear connection
induced by the canonical spray $2G^{i}=\gamma _{jk}^{i}y^{j}y^{k}$ from the
classical theory of Finsler and Lagrange spaces. For more details, please
consult [3], [5].
\end{remark}

Now, let us consider that $\nabla $ is a fixed $\Gamma $-linear connection
on the 1-jet space $E=J^{1}(\mathbb{R},M)$, which is defined by the adapted
local coefficients%
\begin{equation}
\nabla \Gamma =\left( \bar{G}_{11}^{1},\;G_{i1}^{k},\;G_{(1)(j)1}^{(i)(1)},\;%
\bar{L}_{1j}^{1},\;L_{ij}^{k},\;L_{(1)(j)k}^{(i)(1)},\;\bar{C}%
_{1(k)}^{1(1)},\;C_{i(k)}^{j(1)},\;C_{(1)(j)(k)}^{(i)(1)(1)}\right) .
\label{Nabla-Gamma}
\end{equation}

\begin{definition}
A geometrical object $D=\left( D_{1k(1)(l)...}^{1i(j)(1)...}\right) $ on the
1-jet vector bundle $E=J^{1}(\mathbb{R},M)$, whose local components
transform by the rules%
\begin{equation*}
D_{1k(1)(l)...}^{1i(j)(1)...}=\widetilde{D}_{1r(1)(s)...}^{1p(m)(1)...}\frac{%
dt}{d\widetilde{t}}\frac{\partial x^{i}}{\partial \widetilde{x}^{p}}\left( 
\frac{\partial x^{j}}{\partial \widetilde{x}^{m}}\frac{d\widetilde{t}}{dt}%
\right) \frac{d\widetilde{t}}{dt}\frac{\partial \widetilde{x}^{r}}{\partial
x^{k}}\left( \frac{\partial \widetilde{x}^{s}}{\partial x^{l}}\frac{dt}{d%
\widetilde{t}}\right) ...,
\end{equation*}%
is called a \textbf{d-tensor field}.
\end{definition}

\begin{example}
\label{Liouville} The geometrical object $\mathbf{C}=\left( \mathbf{C}%
_{(1)}^{(i)}\right) ,$ where $\mathbf{C}_{(1)}^{(i)}=y_{1}^{i},$ represents
a d-tensor field on the 1-jet space $E=J^{1}(\mathbb{R},M)$. This is called
the \textbf{canonical Liouville d-tensor field} of the 1-jet vector bundle $E
$. Remark that the d-tensor field $\mathbf{C}$ naturally generalizes the
Liouville vector field [5]%
\begin{equation*}
\mathbb{C}=y^{i}\frac{\partial }{\partial y^{i}},
\end{equation*}%
used in the Lagrangian geometry of the tangent bundle $TM$.
\end{example}

The $\Gamma $-linear connection $\nabla $ naturally induces a linear
connection on the set of the d-tensors of the 1-jet vector bundle $E$, in
the following way: $-$ starting with $X\in \mathcal{X}(E)$ a vector field
and $D$ a d-tensor field on $E$, locally expressed by%
\begin{equation*}
\begin{array}{l}
\medskip {X=X^{1}{\dfrac{\delta }{\delta t}}+X^{r}{\dfrac{\delta }{\delta
x^{r}}}+X_{(1)}^{(r)}{\dfrac{\partial }{\partial y_{1}^{r}}},} \\ 
{D=D_{1k(1)(l)\ldots }^{1i(j)(1)\ldots }{\dfrac{\delta }{\delta t}}\otimes {%
\dfrac{\delta }{\delta x^{i}}}\otimes {\dfrac{\partial }{\partial y_{1}^{j}}}%
\otimes dt\otimes dx^{k}\otimes \delta y_{1}^{l}\ldots ,}%
\end{array}%
\end{equation*}%
we introduce the covariant derivative%
\begin{equation*}
\begin{array}{l}
\medskip \nabla _{X}D=X^{1}\nabla _{\dfrac{\delta }{\delta t}}D+X^{p}\nabla
_{\dfrac{\delta }{\delta x^{p}}}D+X_{(1)}^{(p)}\nabla _{\dfrac{\partial }{%
\partial y_{1}^{p}}}D=\left\{ X^{1}D_{1k(1)(l)\ldots /1}^{1i(j)(1)\ldots
}+X^{p}\cdot \right.  \\ 
\left. \cdot D_{1k(1)(l)\ldots |p}^{1i(j)(1)\ldots
}+X_{(1)}^{(p)}D_{1k(1)(l)\ldots }^{1i(j)(1)\ldots }|_{(p)}^{(1)}\right\} {{%
\dfrac{\delta }{\delta t}}\otimes {\dfrac{\delta }{\delta x^{i}}}\otimes {%
\dfrac{\partial }{\partial y_{1}^{j}}}\otimes dt\otimes dx^{k}\otimes \delta
y_{1}^{l}\ldots ,}%
\end{array}%
\end{equation*}%
where

$(h_{\mathbb{R}})\hspace{6mm}\left\{ 
\begin{array}{l}
\medskip {D_{1k(1)(l)\ldots /1}^{1i(j)(1)\ldots }={\dfrac{\delta
D_{1k(1)(l)\ldots }^{1i(j)(1)\ldots }}{\delta t}}+D_{1k(1)(l)\ldots
}^{1i(j)(1)\ldots }\bar{G}_{11}^{1}+} \\ 
\medskip +D_{1k(1)(l)\ldots }^{1r(j)(1)\ldots }G_{r1}^{i}+D_{1k(1)(l)\ldots
}^{1i(r)(1)\ldots }G_{(1)(r)1}^{(j)(1)}+\ldots - \\ 
-D_{1k(1)(l)\ldots }^{1i(j)(1)\ldots }\bar{G}_{11}^{1}-D_{1r(1)(l)\ldots
}^{1i(j)(1)\ldots }G_{k1}^{r}-D_{1k(1)(r)\ldots }^{1i(j)(1)\ldots
}G_{(1)(l)1}^{(r)(1)}-\ldots ,%
\end{array}%
\right. \medskip$

$(h_{M})\hspace{5mm}\left\{ 
\begin{array}{l}
\medskip {D_{1k(1)(l)\ldots |p}^{1i(j)(1)\ldots }={\dfrac{\delta
D_{1k(1)(l)\ldots }^{1i(j)(1)\ldots }}{\delta x^{p}}}+D_{1k(1)(l)\ldots
}^{1i(j)(1)\ldots }\bar{L}_{1p}^{1}+} \\ 
\medskip +D_{1k(1)(l)\ldots }^{1r(j)(1)\ldots }L_{rp}^{i}+D_{1k(1)(l)\ldots
}^{1i(r)(1)\ldots }L_{(1)(r)p}^{(j)(1)}+\ldots - \\ 
-D_{1k(1)(l)\ldots }^{1i(j)(1)\ldots }\bar{L}_{1p}^{1}-D_{1r(1)(l)\ldots
}^{1i(j)(1)\ldots }L_{kp}^{r}-D_{1k(1)(r)\ldots }^{1i(j)(1)\ldots
}L_{(1)(l)p}^{(r)(1)}-\ldots ,%
\end{array}%
\right. \medskip$

$(v)\hspace{8mm}\left\{ 
\begin{array}{l}
\medskip {D_{1k(1)(l)\ldots }^{1i(j)(1)\ldots }|_{(p)}^{(1)}={\dfrac{%
\partial D_{1k(1)(l)\ldots }^{1i(j)(1)\ldots }}{\partial y_{1}^{p}}}%
+D_{1k(1)(l)\ldots }^{1i(j)(1)\ldots }\bar{C}_{1(p)}^{1(1)}+} \\ 
\medskip +D_{1k(1)(l)\ldots }^{1r(j)(1)\ldots
}C_{r(p)}^{i(1)}+D_{1k(1)(l)\ldots }^{1i(r)(1)\ldots
}C_{(1)(r)(p)}^{(j)(1)(1)}+\ldots - \\ 
-D_{1k(1)(l)\ldots }^{1i(j)(1)\ldots }\bar{C}_{1(p)}^{1(1)}-D_{1r(1)(l)%
\ldots }^{1i(j)(1)\ldots }C_{k(p)}^{r(1)}-D_{1k(1)(r)\ldots
}^{1i(j)(1)\ldots }C_{(1)(l)(p)}^{(r)(1)(1)}-\ldots .%
\end{array}%
\right. $

\begin{definition}
The local derivative operators \textquotedblright $_{/1}$\textquotedblright
, \textquotedblright $_{|p}$\textquotedblright\ and \textquotedblright $%
|_{(p)}^{(1)}$\textquotedblright\ are called the $\mathbb{R}$\textbf{%
-horizontal covariant derivative}\textit{, the }$M$\textbf{-horizontal
covariant derivative} and the \textbf{vertical covariant derivative
associated to the }$\Gamma $\textbf{-linear connection }$\nabla \Gamma $.
These \textit{apply to the local components of an arbitrary d-tensor field
on the 1-jet space }$E=J^{1}(\mathbb{R},M)$.
\end{definition}

\begin{remark}
\textbf{(a)} In the particular case of a function $f(t,x^{k},y_{1}^{k})$ on
the 1-jet space $E=J^{1}(\mathbb{R},M)$ the above covariant derivatives
reduce to%
\begin{equation*}
f_{/1}={{\dfrac{\delta f}{\delta t}}={\dfrac{\partial f}{\partial t}}%
-M_{(1)1}^{(k)}{\dfrac{\partial f}{\partial y_{1}^{k}},\quad}}f_{|p}={{%
\dfrac{\delta f}{\delta x^{p}}}={\dfrac{\partial f}{\partial x^{p}}}%
-N_{(1)p}^{(k)}{\dfrac{\partial f}{\partial y_{1}^{k}},\quad}}f|_{(p)}^{(1)}=%
{{\dfrac{\partial f}{\partial y_{1}^{p}}}}.
\end{equation*}

\textbf{(b)} Starting with a d-vector field $D=Y$ on the 1-jet space $%
E=J^{1}(\mathbb{R},M),$ locally expressed by%
\begin{equation*}
Y=Y^{1}{\frac{\delta }{\delta t}}+Y^{i}{\frac{\delta }{\delta x^{i}}}%
+Y_{(1)}^{(i)}{\frac{\partial }{\partial y_{1}^{i}}},
\end{equation*}%
the following expressions of the local covariant derivatives hold
good:\bigskip

$(h_{\mathbb{R}})\hspace*{6mm}\left\{ 
\begin{array}{l}
\medskip {Y_{/1}^{1}={\dfrac{\delta Y^{1}}{\delta t}}+Y^{1}\bar{G}_{11}^{1}}
\\ 
\medskip {Y_{/1}^{i}={\dfrac{\delta Y^{i}}{\delta t}}+Y^{r}G_{r1}^{i}} \\ 
{Y_{(1)/1}^{(i)}={\dfrac{\delta Y_{(1)}^{(i)}}{\delta t}}%
+Y_{(1)}^{(r)}G_{(1)(r)1}^{(i)(1)}},%
\end{array}%
\right. \medskip$

$(h_{M})\hspace*{5mm}\left\{ 
\begin{array}{l}
\medskip {Y_{|p}^{1}={\dfrac{\delta Y^{1}}{\delta x^{p}}}+Y^{1}\bar{L}%
_{1p}^{1}} \\ 
\medskip {Y_{|p}^{i}={\dfrac{\delta Y^{i}}{\delta x^{p}}}+Y^{r}L_{rp}^{i}}
\\ 
{Y_{(1)|p}^{(i)}={\dfrac{\delta Y_{(1)}^{(i)}}{\delta x^{p}}}%
+Y_{(1)}^{(r)}L_{(1)(r)p}^{(i)(1)}},%
\end{array}%
\right. \medskip$

$(v)\hspace*{8mm}\left\{ 
\begin{array}{l}
\medskip {Y^{1}|_{(p)}^{(1)}={\dfrac{\partial Y^{1}}{\partial y_{1}^{p}}}%
+Y^{1}\bar{C}_{1(p)}^{1(1)}} \\ 
\medskip {Y^{i}|_{(p)}^{(1)}={\dfrac{\partial Y^{i}}{\partial y_{1}^{p}}}%
+Y^{r}C_{r(p)}^{i(1)}} \\ 
{Y_{(1)}^{(i)}|_{(p)}^{(1)}={\dfrac{\partial Y_{(1)}^{(i)}}{\partial
y_{1}^{p}}}+Y_{(1)}^{(r)}C_{(1)(r)(p)}^{(i)(1)(1)}}.%
\end{array}%
\right. \medskip$
\end{remark}

Denoting generically by \textquotedblright $_{:A}$\textquotedblright\ one of
the local covariant derivatives \textquotedblright $_{/1}$\textquotedblright
, \textquotedblright $_{|p}$\textquotedblright\ or \textquotedblright $%
|_{(p)}^{(1)}$\textquotedblright , we obtain the following properties of the
covariant derivative operators:

\begin{proposition}
If$\ T_{...}^{...}$ and $S_{...}^{...}$ are two arbitrary d-tensors on $%
E=J^{1}(\mathbb{R},M),$ then the following statements hold good:\medskip

\textbf{(i)} The local coefficients $T_{\ldots :A}^{\ldots }$ represent the
components of a new d-tensor field on the 1-jet space $E=J^{1}(\mathbb{R},M)$%
.\medskip

\textbf{(ii)} $(T_{\ldots }^{\ldots }+S_{\ldots }^{\ldots })_{:A}=T_{\ldots
:A}^{\ldots }+S_{\ldots :A}^{\ldots }$.\medskip

\textbf{(iii)} $(T_{\ldots }^{\ldots }\otimes S_{\ldots }^{\ldots
})_{:A}=T_{\ldots :A}^{\ldots }\otimes S_{\ldots }^{\ldots }+T_{\ldots
}^{\ldots }\otimes S_{\ldots :A}^{\ldots }$.
\end{proposition}

\section{Torsion and curvature d-tensors}

\hspace{5mm}In the sequel, we will study the torsion tensor $\mathbf{T}:%
\mathcal{X}(E)\times \mathcal{X}(E)\rightarrow \mathcal{X}(E)$ associated to
the $\Gamma $-linear connection $\nabla $, which is given by the formula%
\begin{equation*}
\mathbf{T}(X,Y)=\nabla _{X}Y-\nabla _{Y}X-[X,Y],\quad \forall \;X,Y\in 
\mathcal{X}(E).
\end{equation*}

In order to obtain an adapted local characterization of the torsion tensor $%
\mathbf{T}$ of the $\Gamma $-linear connection $\nabla $, we firstly deduce,
by direct computations, the following simple and important result:

\begin{proposition}
The following identities of the Poisson brackets are true: 
\begin{equation*}
\begin{array}{ll}
\medskip {\left[ {\dfrac{\delta }{\delta t}},{\dfrac{\delta }{\delta t}}%
\right] =0}, & {\left[ {\dfrac{\delta }{\delta t}},{\dfrac{\delta }{\delta
x^{j}}}\right] =R_{(1)1j}^{(r)}{\dfrac{\partial }{\partial y_{1}^{r}}}}, \\ 
\medskip {\left[ {\dfrac{\delta }{\delta t}},{\dfrac{\partial }{\partial
y_{1}^{j}}}\right] ={\dfrac{\partial M_{(1)1}^{(r)}}{\partial y_{1}^{j}}}{%
\dfrac{\partial }{\partial y_{1}^{r}}}}, & {\left[ {\dfrac{\delta }{\delta
x^{i}}},{\dfrac{\delta }{\delta x^{j}}}\right] =R_{(1)ij}^{(r)}{\dfrac{%
\partial }{\partial y_{1}^{r}}}}, \\ 
{\left[ {\dfrac{\delta }{\delta x^{i}}},{\dfrac{\partial }{\partial y_{1}^{j}%
}}\right] ={\dfrac{\partial N_{(1)i}^{(r)}}{\partial y_{1}^{j}}}{\dfrac{%
\partial }{\partial y_{1}^{r}}}}, & {\left[ {\dfrac{\partial }{\partial
y_{1}^{i}}},{\dfrac{\partial }{\partial y_{1}^{j}}}\right] =0},%
\end{array}%
\end{equation*}%
\medskip where $M_{(1)1}^{(r)}$ and $N_{(1)i}^{(r)}$ are the local
coefficients of the nonlinear connection $\Gamma ,$ while the components $%
R_{(1)1j}^{(r)}\;$and$\;R_{(1)ij}^{(r)}$ are d-tensors given by the formulas%
\begin{equation}
\begin{array}{l}
R_{(1)1j}^{(r)}={\dfrac{\delta M_{(1)1}^{(r)}}{\delta x^{j}}}-{\dfrac{\delta
N_{(1)j}^{(r)}}{\delta t}},\medskip \\ 
R_{(1)ij}^{(r)}={\dfrac{\delta N_{(1)i}^{(r)}}{\delta x^{j}}}-{\dfrac{\delta
N_{(1)j}^{(r)}}{\delta x^{i}}}.%
\end{array}
\label{Poisson brackets jet}
\end{equation}
\end{proposition}

In these conditions, working with a basis of vector fields, adapted to the
nonlinear connection 
\begin{equation*}
\Gamma =\left( M_{(1)1}^{(i)},N_{(1)j}^{(i)}\right)
\end{equation*}%
on the 1-jet space $E=J^{1}(\mathbb{R},M)$, by local computations, we obtain

\begin{theorem}
The torsion tensor $\mathbf{T}$ of the $\Gamma $-linear connection (\ref%
{Nabla-Gamma}) is determined by the following adapted torsion d-tensors:%
\begin{equation*}
{h_{\mathbb{R}}\mathbf{T}\left( {\frac{\delta }{\delta t}},{\frac{\delta }{%
\delta t}}\right) =0,\quad h_{M}\mathbf{T}\left( {\frac{\delta }{\delta t}},{%
\frac{\delta }{\delta t}}\right) =0,\quad v\mathbf{T}\left( {\frac{\delta }{%
\delta t}},{\frac{\delta }{\delta t}}\right) =0},
\end{equation*}%
\begin{equation*}
{h_{\mathbb{R}}\mathbf{T}\left( {\frac{\delta }{\delta x^{j}}},{\frac{\delta 
}{\delta t}}\right) =\bar{T}_{1j}^{1}{\frac{\delta }{\delta t}},\quad h_{M}%
\mathbf{T}\left( {\frac{\delta }{\delta x^{j}}},{\frac{\delta }{\delta t}}%
\right) =T_{1j}^{r}{\frac{\delta }{\delta x^{r}}},}
\end{equation*}%
\begin{equation*}
{v\mathbf{T}\left( {\frac{\delta }{\delta x^{j}}},{\frac{\delta }{\delta t}}%
\right) =R_{(1)1j}^{(r)}{\frac{\partial }{\partial y_{1}^{r}}}},
\end{equation*}%
\begin{equation*}
{h_{\mathbb{R}}\mathbf{T}\left( {\frac{\delta }{\delta x^{j}}},{\frac{\delta 
}{\delta x^{i}}}\right) =0},\quad {h_{M}\mathbf{T}\left( {\frac{\delta }{%
\delta x^{j}}},{\frac{\delta }{\delta x^{i}}}\right) =T_{ij}^{r}{\frac{%
\delta }{\delta x^{r}}}},
\end{equation*}%
\begin{equation*}
{v\mathbf{T}\left( {\frac{\delta }{\delta x^{j}}},{\frac{\delta }{\delta
x^{i}}}\right) =R_{(1)ij}^{(r)}{\frac{\partial }{\partial y_{1}^{r}}}},
\end{equation*}%
\begin{equation*}
{h_{\mathbb{R}}\mathbf{T}\left( {\frac{\partial }{\partial y_{1}^{j}}},{%
\frac{\delta }{\delta t}}\right) =\bar{P}_{1(j)}^{1(1)}{\frac{\delta }{%
\delta t}}},\quad {h_{M}\mathbf{T}\left( {\frac{\partial }{\partial y_{1}^{j}%
}},{\frac{\delta }{\delta t}}\right) =0},
\end{equation*}%
\begin{equation*}
{v\mathbf{T}\left( {\frac{\partial }{\partial y_{1}^{j}}},{\frac{\delta }{%
\delta t}}\right) =P_{(1)1(j)}^{(r)\text{ \ }(1)}{\frac{\partial }{\partial
y_{1}^{r}}}},
\end{equation*}%
\begin{equation*}
{h_{\mathbb{R}}\mathbf{T}\left( {\frac{\partial }{\partial y_{1}^{j}}},{%
\frac{\delta }{\delta x^{i}}}\right) =0},\quad {h_{M}\mathbf{T}\left( {\frac{%
\partial }{\partial y_{1}^{j}}},{\frac{\delta }{\delta x^{i}}}\right)
=P_{i(j)}^{r(1)}{\frac{\delta }{\delta x^{r}}}},
\end{equation*}%
\begin{equation*}
{v\mathbf{T}\left( {\frac{\partial }{\partial y_{1}^{j}}},{\frac{\delta }{%
\delta x^{i}}}\right) =P_{(1)i(j)}^{(r)\;(1)}{\frac{\partial }{\partial
y_{1}^{r}}}},
\end{equation*}%
\begin{equation*}
{h_{\mathbb{R}}\mathbf{T}\left( {\frac{\partial }{\partial y_{1}^{j}}},{%
\frac{\partial }{\partial y_{1}^{i}}}\right) =0},\quad {h_{M}\mathbf{T}%
\left( {\frac{\partial }{\partial y_{1}^{j}}},{\frac{\partial }{\partial
y_{1}^{i}}}\right) =0},
\end{equation*}%
\begin{equation*}
{v\mathbf{T}\left( {\frac{\partial }{\partial y_{1}^{j}}},{\frac{\partial }{%
\partial y_{1}^{i}}}\right) =S_{(1)(i)(j)}^{(r)(1)(1)}{\frac{\partial }{%
\partial y_{1}^{r}}}},
\end{equation*}%
where 
\begin{equation*}
\begin{array}{cccc}
\bar{T}_{1j}^{1}=\bar{L}_{1j}^{1}, & T_{1j}^{r}=-G_{j1}^{r}, & 
T_{ij}^{r}=L_{ij}^{r}-L_{ji}^{r}, & \bar{P}_{1(j)}^{1(1)}=\bar{C}%
_{1(j)}^{1(1)},%
\end{array}%
\end{equation*}%
\medskip 
\begin{equation*}
\begin{array}{cc}
P_{i(j)}^{r(1)}=C_{i(j)}^{r(1)}, & 
S_{(1)(i)(j)}^{(r)(1)(1)}=C_{(1)(i)(j)}^{(r)(1)(1)}-C_{(1)(j)(i)}^{(r)(1)(1)},%
\end{array}%
\end{equation*}%
\medskip 
\begin{equation*}
\begin{array}{cc}
{P_{(1)1(j)}^{(r)\text{ \ }(1)}={\dfrac{\partial M_{(1)1}^{(r)}}{\partial
y_{1}^{j}}}-G_{(1)(j)1}^{(r)(1)}}, & {P_{(1)i(j)}^{(r)\text{ }(1)}={\dfrac{%
\partial N_{(1)i}^{(r)}}{\partial y_{1}^{j}}}-L_{(1)(j)i}^{(r)(1)}},\medskip%
\end{array}%
\end{equation*}%
\medskip and the d-tensors $R_{(1)1j}^{(r)}$ and $R_{(1)ij}^{(r)}$ are given
by (\ref{Poisson brackets jet}).
\end{theorem}

\begin{corollary}
The torsion tensor $\mathbf{T}$ of an arbitrary $\Gamma $-linear connection $%
\nabla $ on the 1-jet space $E=J^{1}(\mathbb{R},M)$ is determined by \textbf{%
ten} effective adapted local torsion d-tensors, which we arrange in the
following table:%
\begin{equation*}
\begin{tabular}{|c|c|c|c|}
\hline
& $h_{\mathbb{R}}$ & $h_{M}$ & $v$ \\ \hline
$h_{\mathbb{R}}h_{\mathbb{R}}$ & $0$ & $0$ & $0$ \\ \hline
$h_{M}h_{\mathbb{R}}$ & $\bar{T}_{1j}^{1}$ & $T_{1j}^{r}$ & $R_{(1)1j}^{(r)}$
\\ \hline
$h_{M}h_{M}$ & $0$ & $T_{ij}^{r}$ & $R_{(1)ij}^{(r)}$ \\ \hline
$vh_{\mathbb{R}}$ & $\bar{P}_{1(j)}^{1(1)}$ & $0$ & ${P_{(1)1(j)}^{(r)\text{
\ }(1)}}$ \\ \hline
$vh_{M}$ & $0$ & $P_{i(j)}^{r(1)}$ & ${P_{(1)i(j)}^{(r)\;(1)}}$ \\ \hline
$vv$ & $0$ & $0$ & $S_{(1)(i)(j)}^{(r)(1)(1)}$ \\ \hline
\end{tabular}%
\end{equation*}
\end{corollary}

\begin{example}
In the particular case of the Berwald $\mathring{\Gamma}$-linear connection $%
B\mathring{\Gamma},$ associated to the semi-Riemannian metrics $h_{11}(t)$
and $\varphi _{ij}(x),$ all torsion d-tensors vanish, except%
\begin{equation*}
R_{(1)ij}^{(k)}=\mathfrak{R}_{mij}^{k}y_{1}^{m},
\end{equation*}%
where $\mathfrak{R}_{mij}^{k}(x)$ are the classical local curvature tensors
of the spatial semi-Rie\-ma\-nni\-an metric $\varphi _{ij}(x)$.
\end{example}

In order to study the curvature of the $\Gamma $-linear connection $\nabla $%
, we recall that the curvature tensor $\mathbf{R}$ of $\nabla $ is given by
the formula%
\begin{equation*}
\mathbf{R}(X,Y)Z=\nabla _{X}\nabla _{Y}Z-\nabla _{Y}\nabla _{X}Z-\nabla
_{\lbrack X,Y]}Z,\quad \forall \;X,Y,Z\in \mathcal{X}(E).
\end{equation*}

Using again a basis of vector fields adapted to the nonlinear connection $%
\Gamma $, together with the properties of the $\Gamma $-linear connection $%
\nabla $, by direct computations, we obtain

\begin{theorem}
The curvature tensor $\mathbf{R}$ associated to the $\Gamma $-linear
connection (\ref{Nabla-Gamma}) is determined by \textbf{fifteen} effective
adapted local curvature d-tensors%
\begin{equation*}
\mathbf{R}{\left( {\frac{\delta }{\delta t}},{\frac{\delta }{\delta t}}%
\right) {\frac{\delta }{\delta t}}=0},\quad \mathbf{R}{\left( {\frac{\delta 
}{\delta t}},{\frac{\delta }{\delta t}}\right) {\frac{\delta }{\delta x^{i}}}%
=0},\quad \mathbf{R}{\left( {\frac{\delta }{\delta t}},{\frac{\delta }{%
\delta t}}\right) {\frac{\partial }{\partial y_{1}^{i}}}=0},
\end{equation*}%
\begin{equation*}
\mathbf{R}{\left( {\frac{\delta }{\delta x^{k}}},{\frac{\delta }{\delta t}}%
\right) {\frac{\delta }{\delta t}}=\bar{R}_{11k}^{1}{\frac{\delta }{\delta t}%
}},\quad \mathbf{R}{\left( {\frac{\delta }{\delta x^{k}}},{\frac{\delta }{%
\delta t}}\right) {\frac{\delta }{\delta x^{i}}}=R_{i1k}^{l}{\frac{\delta }{%
\delta x^{l}}}},
\end{equation*}%
\begin{equation*}
\mathbf{R}{\left( {\frac{\delta }{\delta x^{k}}},{\frac{\delta }{\delta t}}%
\right) {\frac{\partial }{\partial y_{1}^{i}}}=R_{(1)(i)1k}^{(l)(1)}{\frac{%
\partial }{\partial y_{1}^{l}}}},
\end{equation*}%
\begin{equation*}
\mathbf{R}{\left( {\frac{\delta }{\delta x^{k}}},{\frac{\delta }{\delta x^{j}%
}}\right) {\frac{\delta }{\delta t}}=\bar{R}_{1jk}^{1}{\frac{\delta }{\delta
t}}},\quad \mathbf{R}{\left( {\frac{\delta }{\delta x^{k}}},{\frac{\delta }{%
\delta x^{j}}}\right) {\frac{\delta }{\delta x^{i}}}=R_{ijk}^{l}{\frac{%
\delta }{\delta x^{l}}}},
\end{equation*}%
\begin{equation*}
\mathbf{R}{\left( {\frac{\delta }{\delta x^{k}}},{\frac{\delta }{\delta x^{j}%
}}\right) {\frac{\partial }{\partial y_{1}^{i}}}=R_{(1)(i)jk}^{(l)(1)}{\frac{%
\partial }{\partial y_{1}^{l}}}},
\end{equation*}%
\begin{equation*}
\mathbf{R}{\left( {\frac{\partial }{\partial y_{1}^{k}}},{\frac{\delta }{%
\delta t}}\right) {\frac{\delta }{\delta t}}=\bar{P}_{11(k)}^{1\;\;(1)}{%
\frac{\delta }{\delta t}}},\quad \mathbf{R}{\left( {\frac{\partial }{%
\partial y_{1}^{k}}},{\frac{\delta }{\delta t}}\right) {\frac{\delta }{%
\delta x^{i}}}=P_{i1(k)}^{l\;\;(1)}{\frac{\delta }{\delta x^{l}}}},
\end{equation*}%
\begin{equation*}
\mathbf{R}{\left( {\frac{\partial }{\partial y_{1}^{k}}},{\frac{\delta }{%
\delta t}}\right) {\frac{\partial }{\partial y_{1}^{i}}}%
=P_{(1)(i)1(k)}^{(l)(1)\;(1)}{\frac{\partial }{\partial y_{1}^{l}}}},
\end{equation*}%
\begin{equation*}
\mathbf{R}{\left( {\frac{\partial }{\partial y_{1}^{k}}},{\frac{\delta }{%
\delta x^{j}}}\right) {\frac{\delta }{\delta t}}=\bar{P}_{1j(k)}^{1\;\;(1)}{%
\frac{\delta }{\delta t}}},\quad \mathbf{R}{\left( {\frac{\partial }{%
\partial y_{1}^{k}}},{\frac{\delta }{\delta x^{j}}}\right) {\frac{\delta }{%
\delta x^{i}}}=P_{ij(k)}^{l\;\;(1)}{\frac{\delta }{\delta x^{l}}}},
\end{equation*}%
\begin{equation*}
\mathbf{R}{\left( {\frac{\partial }{\partial y_{1}^{k}}},{\frac{\delta }{%
\delta x^{j}}}\right) {\frac{\partial }{\partial y_{1}^{i}}}%
=P_{(1)(i)j(k)}^{(l)(1)\;(1)}{\frac{\partial }{\partial y_{1}^{l}}}},
\end{equation*}%
\begin{equation*}
\mathbf{R}{\left( {\frac{\partial }{\partial y_{1}^{k}}},{\frac{\partial }{%
\partial y_{1}^{j}}}\right) {\frac{\delta }{\delta t}}=\bar{S}%
_{1(j)(k)}^{1(1)(1)}{\frac{\delta }{\delta t}}},\quad \mathbf{R}{\left( {%
\frac{\partial }{\partial y_{1}^{k}}},{\frac{\partial }{\partial y_{1}^{j}}}%
\right) {\frac{\delta }{\delta x^{i}}}=S_{i(j)(k)}^{l(1)(1)}{\frac{\delta }{%
\delta x^{l}}}},
\end{equation*}%
\begin{equation*}
\mathbf{R}{\left( {\frac{\partial }{\partial y_{1}^{k}}},{\frac{\partial }{%
\partial y_{1}^{j}}}\right) {\frac{\partial }{\partial y_{1}^{i}}}%
=S_{(1)(i)(j)(k)}^{(l)(1)(1)(1)}{\frac{\partial }{\partial y_{1}^{l}}}},
\end{equation*}%
whose local components we arrange in the following table:%
\begin{equation*}
\begin{tabular}{|c|c|c|c|}
\hline
& $h_{\mathbb{R}}$ & $h_{M}$ & $v$ \\ \hline
$h_{\mathbb{R}}h_{\mathbb{R}}$ & $0$ & $0$ & $0$ \\ \hline
$h_{M}h_{\mathbb{R}}$ & $\bar{R}_{11k}^{1}$ & $R_{i1k}^{l}$ & $%
R_{(1)(i)1k}^{(l)(1)}$ \\ \hline
$h_{M}h_{M}$ & $\bar{R}_{1jk}^{1}$ & $R_{ijk}^{l}$ & $R_{(1)(i)jk}^{(l)(1)}$
\\ \hline
$vh_{\mathbb{R}}$ & $\bar{P}_{11(k)}^{1\;\;(1)}$ & $P_{i1(k)}^{l\;\;(1)}$ & $%
P_{(1)(i)1(k)}^{(l)(1)\;(1)}$ \\ \hline
$vh_{M}$ & $\bar{P}_{1j(k)}^{1\;\;(1)}$ & $P_{ij(k)}^{l\;\;(1)}$ & $%
P_{(1)(i)j(k)}^{(l)(1)\;(1)}$ \\ \hline
$vv$ & $\bar{S}_{1(j)(k)}^{1(1)(1)}$ & $S_{i(j)(k)}^{l(1)(1)}$ & $%
S_{(1)(i)(j)(k)}^{(l)(1)(1)(1)}$ \\ \hline
\end{tabular}%
\end{equation*}
\end{theorem}

Moreover, by a laborious local computations, we deduce the following result:

\begin{theorem}
The expressions of the preceding local curvature d-tensors are given by:

\begin{itemize}
\item $h_{\mathbb{R}}$-components\medskip

$\ \mathbf{1.}$ ${\bar{R}_{11k}^{1}={\dfrac{\delta \bar{G}_{11}^{1}}{\delta
x^{k}}}-{\dfrac{\delta \bar{L}_{1k}^{1}}{\delta t}}+\bar{C}%
_{1(r)}^{1(1)}R_{(1)1k}^{(r)}}\medskip $

$\ \mathbf{2.}$ ${\bar{R}_{1jk}^{1}={\dfrac{\delta \bar{L}_{1j}^{1}}{\delta
x^{k}}}-{\dfrac{\delta \bar{L}_{1k}^{1}}{\delta x^{j}}}+\bar{C}%
_{1(r)}^{1(1)}R_{(1)jk}^{(r)}}\medskip $

$\ \mathbf{3.}$ ${\bar{P}_{11(k)}^{1\;\;(1)}={\dfrac{\partial \bar{G}%
_{11}^{1}}{\partial y_{1}^{k}}}-\bar{C}_{1(k)/1}^{1(1)}+\bar{C}%
_{1(r)}^{1(1)}P_{(1)1(k)}^{(r)\;\;(1)}}\medskip $

$\ \mathbf{4.}$ ${\bar{P}_{1j(k)}^{1\;\;(1)}={\dfrac{\partial \bar{L}%
_{1j}^{1}}{\partial y_{1}^{k}}}-\bar{C}_{1(k)|j}^{1(1)}+\bar{C}%
_{1(r)}^{1(1)}P_{(1)j(k)}^{(r)\;\;(1)}}\medskip $

$\ \mathbf{5.}$ ${\bar{S}_{1(j)(k)}^{1(1)(1)}={\dfrac{\partial \bar{C}%
_{1(j)}^{1(1)}}{\partial y_{1}^{k}}}-{\dfrac{\partial \bar{C}_{1(k)}^{1(1)}}{%
\partial y_{1}^{j}}}}\medskip $

\item $h_{M}$-components\medskip

$\ \mathbf{6.}$ ${R_{i1k}^{l}={\dfrac{\delta G_{i1}^{l}}{\delta x^{k}}}-{%
\dfrac{\delta L_{ik}^{l}}{\delta t}}%
+G_{i1}^{r}L_{rk}^{l}-L_{ik}^{r}G_{r1}^{l}+C_{i(r)}^{l(1)}R_{(1)1k}^{(r)}}%
\medskip $

$\ \mathbf{7.}$ ${R_{ijk}^{l}={\dfrac{\delta L_{ij}^{l}}{\delta x^{k}}}-{%
\dfrac{\delta L_{ik}^{l}}{\delta x^{j}}}%
+L_{ij}^{r}L_{rk}^{l}-L_{ik}^{r}L_{rj}^{l}+C_{i(r)}^{l(1)}R_{(1)jk}^{(r)}}%
\medskip $

\textbf{\ }$\mathbf{8.}$ ${P_{i1(k)}^{l\;\;(1)}={\dfrac{\partial G_{i1}^{l}}{%
\partial y_{1}^{k}}}-C_{i(k)/1}^{l(1)}+C_{i(r)}^{l(1)}P_{(1)1(k)}^{(r)\;%
\;(1)}}\medskip $

\textbf{\ }$\mathbf{9.}$ ${P_{ij(k)}^{l\;\;(1)}={\dfrac{\partial L_{ij}^{l}}{%
\partial y_{1}^{k}}}-C_{i(k)|j}^{l(1)}+C_{i(r)}^{l(1)}P_{(1)j(k)}^{(r)\;%
\;(1)}}\medskip $

$\mathbf{10.}$ ${S_{i(j)(k)}^{l(1)(1)}={\dfrac{\partial C_{i(j)}^{l(1)}}{%
\partial y_{1}^{k}}}-{\dfrac{\partial C_{i(k)}^{l(1)}}{\partial y_{1}^{j}}}%
+C_{i(j)}^{r(1)}C_{r(k)}^{l(1)}-C_{i(k)}^{r(1)}C_{r(j)}^{l(1)}}\medskip $

\item $v$-components\medskip

$\mathbf{11.}$ ${R_{(1)(i)1k}^{(l)(1)}={\dfrac{\delta G_{(1)(i)1}^{(l)(1)}}{%
\delta x^{k}}}-{\dfrac{\delta L_{(1)(i)k}^{(l)(1)}}{\delta t}}%
+G_{(1)(i)1}^{(r)(1)}L_{(1)(r)k}^{(l)(1)}-}%
L_{(1)(i)k}^{(r)(1)}G_{(1)(r)1}^{(l)(1)}+\medskip $

$\mbox{\hspace{24mm}}+C_{(1)(i)(r)}^{(l)(1)(1)}R_{(1)1k}^{(r)}\medskip $

$\mathbf{12.}$ ${R_{(1)(i)jk}^{(l)(1)}={\dfrac{\delta L_{(1)(i)j}^{(l)(1)}}{%
\delta x^{k}}}-{\dfrac{\delta L_{(1)(i)k}^{(l)(1)}}{\delta x^{j}}}%
+L_{(1)(i)j}^{(r)(1)}L_{(1)(r)k}^{(l)(1)}-}%
L_{(1)(i)k}^{(r)(1)}L_{(1)(r)j}^{(l)(1)}+\medskip $

$\mbox{\hspace{24mm}}+C_{(1)(i)(r)}^{(l)(1)(1)}R_{(1)jk}^{(r)}\medskip $

$\mathbf{13.}$ ${P_{(1)(i)1(k)}^{(l)(1)\;\;(1)}={\dfrac{\partial
G_{(1)(i)1}^{(l)(1)}}{\partial y_{1}^{k}}}%
-C_{(1)(i)(k)/1}^{(l)(1)(1)}+C_{(1)(i)(r)}^{(l)(1)(1)}P_{(1)1(k)}^{(r)\;%
\;(1)}}\medskip $

$\mathbf{14.}$ ${P_{(1)(i)j(k)}^{(l)(1)\;\;(1)}={\dfrac{\partial
L_{(1)(i)j}^{(l)(1)}}{\partial y_{1}^{k}}}%
-C_{(1)(i)(k)|j}^{(l)(1)(1)}+C_{(1)(i)(r)}^{(l)(1)(1)}P_{(1)j(k)}^{(r)\;%
\;(1)}}\medskip $

$\mathbf{15.}$ ${S_{(1)(i)(j)(k)}^{(l)(1)(1)(1)}={\dfrac{\partial
C_{(1)(i)(j)}^{(l)(1)(1)}}{\partial y_{1}^{k}}}-{\dfrac{\partial
C_{(1)(i)(k)}^{(l)(1)(1)}}{\partial y_{1}^{j}}}%
+C_{(1)(i)(j)}^{(r)(1)(1)}C_{(1)(r)(k)}^{(l)(1)(1)}-}\medskip $

$\mbox{\hspace{29mm}}-C_{(1)(i)(k)}^{(r)(1)(1)}C_{(1)(r)(j)}^{(l)(1)(1)}.$
\end{itemize}
\end{theorem}

\begin{example}
In the case of the Berwald $\mathring{\Gamma}$-linear connection $B\mathring{%
\Gamma},$ associated to the pair of semi-Riemannian metrics $%
(h_{ab}(t),\varphi _{ij}(x)),$ all local curvature d-tensors vanish, except%
\begin{equation*}
R_{ijk}^{l}=\mathfrak{R}_{ijk}^{l},\quad {R_{(1)(i)jk}^{(l)(1)}=}\mathfrak{R}%
_{ijk}^{l},
\end{equation*}%
where $\mathfrak{R}_{ijk}^{l}(x)$ are the classical local curvature tensors
of the spatial semi-Rie\-ma\-nni\-an metric $\varphi _{ij}(x)$.
\end{example}

\section{Ricci identities and deflection d-tensors}

\hspace{5mm}Using the properties of a $\Gamma $-linear connection $\nabla $
given by (\ref{Nabla-Gamma}), together with the definitions of its torsion
tensor $\mathbf{T}$ and its curvature tensor $\mathbf{R}$, we can prove the
following important result which is used in the Lagrangian geometrical
theory of the relativistic time dependent electromagnetism, in order to
describe its generalized Maxwell equations. For more details, please consult
[8].

\begin{theorem}
\label{ricci} If $X$ is an arbitrary d-vector field on the 1-jet vector
bundle $E=J^{1}(\mathbb{R},M),$ locally expressed by%
\begin{equation*}
{X=X^{1}{\frac{\delta }{\delta t}}+X^{i}{\frac{\delta }{\delta x^{i}}}%
+X_{(1)}^{(i)}{\frac{\partial }{\partial y_{1}^{i}},}}
\end{equation*}%
then the following \textbf{Ricci identities} of the $\Gamma $-linear
connection $\nabla $ are true:\medskip

$(h_{\mathbb{R}})$ $\ \hspace{1mm}\left\{ 
\begin{array}{l}
\medskip X_{/1|k}^{1}-X_{|k/1}^{1}=X^{1}\bar{R}_{11k}^{1}-X_{/1}^{1}\bar{T}%
_{1k}^{1}-X_{|r}^{1}T_{1k}^{r}-X^{1}|_{(r)}^{(1)}R_{(1)1k}^{(r)} \\ 
\medskip X_{|j|k}^{1}-X_{|k|j}^{1}=X^{1}\bar{R}%
_{1jk}^{1}-X_{|r}^{1}T_{jk}^{r}-X^{1}|_{(r)}^{(1)}R_{(1)jk}^{(r)} \\ 
\medskip X_{/1}^{1}|_{(k)}^{(1)}-X^{1}|_{(k)/1}^{(1)}=X^{1}\bar{P}%
_{11(k)}^{1\;\;(1)}-X_{/1}^{1}\bar{C}%
_{1(k)}^{1(1)}-X^{1}|_{(r)}^{(1)}P_{(1)1(k)}^{(r)\;\;(1)} \\ 
\medskip X_{|j}^{1}|_{(k)}^{(1)}-X^{1}|_{(k)|j}^{(1)}=X^{1}\bar{P}%
_{1j(k)}^{1\;%
\;(1)}-X_{|r}^{1}C_{j(k)}^{r(1)}-X^{1}|_{(r)}^{(1)}P_{(1)j(k)}^{(r)\;\;(1)}
\\ 
\medskip
X^{1}|_{(j)}^{(1)}|_{(k)}^{(1)}-X^{1}|_{(k)}^{(1)}|_{(j)}^{(1)}=X^{1}\bar{S}%
_{1(j)(k)}^{1(1)(1)}-X^{1}|_{(r)}^{(1)}S_{(1)(j)(k)}^{(r)(1)(1)}%
\end{array}%
\right. \medskip $

$(h_{M})$ $\ \left\{ 
\begin{array}{l}
\medskip X_{/1|k}^{i}-X_{|k/1}^{i}=X^{r}R_{r1k}^{i}-X_{/1}^{i}\bar{T}%
_{1k}^{1}-X_{|r}^{i}T_{1k}^{r}-X^{i}|_{(r)}^{(1)}R_{(1)1k}^{(r)} \\ 
\medskip
X_{|j|k}^{i}-X_{|k|j}^{i}=X^{r}R_{rjk}^{i}-X_{|r}^{i}T_{jk}^{r}-X^{i}|_{(r)}^{(1)}R_{(1)jk}^{(r)}
\\ 
\medskip
X_{/1}^{i}|_{(k)}^{(1)}-X^{i}|_{(k)/1}^{(1)}=X^{r}P_{r1(k)}^{i\;%
\;(1)}-X_{/1}^{i}\bar{C}_{1(k)}^{1(1)}-X^{i}|_{(r)}^{(1)}P_{(1)1(k)}^{(r)\;%
\;(1)} \\ 
\medskip
X_{|j}^{i}|_{(k)}^{(1)}-X^{i}|_{(k)|j}^{(1)}=X^{r}P_{rj(k)}^{i\;%
\;(1)}-X_{|r}^{i}C_{j(k)}^{r(1)}-X^{i}|_{(r)}^{(1)}P_{(1)j(k)}^{(r)\;\;(1)}
\\ 
\medskip
X^{i}|_{(j)}^{(1)}|_{(k)}^{(1)}-X^{i}|_{(k)}^{(1)}|_{(j)}^{(1)}=X^{r}S_{r(j)(k)}^{i(1)(1)}-X^{i}|_{(r)}^{(1)}S_{(1)(j)(k)}^{(r)(1)(1)}%
\end{array}%
\right. \medskip $

$(v)$ $\ \hspace{4mm}\left\{ 
\begin{array}{l}
\medskip
X_{(1)/1|k}^{(i)}-X_{(1)|k/1}^{(i)}=X_{(1)}^{(r)}R_{(1)(r)1k}^{(i)(1)}-X_{(1)/1}^{(i)}%
\bar{T}_{1k}^{1}-X_{(1)|r}^{(i)}T_{1k}^{r}- \\ 
\medskip \mbox{\hspace{30mm}}-X_{(1)}^{(i)}|_{(r)}^{(1)}R_{(1)1k}^{(r)} \\ 
\medskip
X_{(1)|j|k}^{(i)}-X_{(1)|k|j}^{(i)}=X_{(1)}^{(r)}R_{(1)(r)jk}^{(i)(1)}-X_{(1)|r}^{(i)}T_{jk}^{r}-X_{(1)}^{(i)}|_{(r)}^{(1)}R_{(1)jk}^{(r)}
\\ 
\medskip
X_{(1)/1}^{(i)}|_{(k)}^{(1)}-X_{(1)}^{(i)}|_{(k)/1}^{(1)}=X_{(1)}^{(r)}P_{(1)(r)1(k)}^{(i)(1)\;\;(1)}-X_{(1)/1}^{(i)}%
\bar{C}_{1(k)}^{1(1)}- \\ 
\medskip \mbox{\hspace{35mm}}-X_{(1)}^{(i)}|_{(r)}^{(1)}P_{(1)1(k)}^{(r)\;%
\;(1)} \\ 
\medskip
X_{(1)|j}^{(i)}|_{(k)}^{(1)}-X_{(1)}^{(i)}|_{(k)|j}^{(1)}=X_{(1)}^{(r)}P_{(1)(r)j(k)}^{(i)(1)\;\;(1)}-X_{(1)|r}^{(i)}C_{j(k)}^{r(1)}-
\\ 
\medskip \mbox{\hspace{33mm}}-X_{(1)}^{(i)}|_{(r)}^{(1)}P_{(1)j(k)}^{(r)\;%
\;(1)} \\ 
\medskip
X_{(1)}^{(i)}|_{(j)}^{(1)}|_{(k)}^{(1)}-X_{(1)}^{(i)}|_{(k)}^{(1)}|_{(j)}^{(1)}=X_{(1)}^{(r)}S_{(1)(r)(j)(k)}^{(i)(1)(1)(1)}-X_{(1)}^{(i)}|_{(r)}^{(1)}S_{(1)(j)(k)}^{(r)(1)(1)}.%
\end{array}%
\right. $
\end{theorem}

\begin{proof}
Let $(Y_{A})$ and $(\omega ^{A})$, where $A\in \left\{ 1,i,{\QATOP{(1)}{(i)}}%
\right\} $, be the dual bases adapted to the nonlinear connection $\Gamma $
and let $X=X^{F}Y_{F}$ be a distinguished vector field on the 1-jet space $%
E=J^{1}(\mathbb{R},M)$. In this context, using the equalities%
\begin{equation*}
\begin{array}{ll}
\mathbf{1.} & \nabla _{Y_{C}}Y_{B}=\Gamma _{BC}^{F}Y_{F},\qquad\mathbf{2.}%
\text{\hspace{3mm}}[Y_{B},Y_{C}]=R_{BC}^{F}Y_{F},\medskip \\ 
\mathbf{3.} & \mathbf{T}(Y_{C},Y_{B})=\mathbf{T}_{BC}^{F}Y_{F}=\{\Gamma
_{BC}^{F}-\Gamma _{CB}^{F}-R_{CB}^{F}\}Y_{F},\medskip \\ 
\mathbf{4.} & \mathbf{R}(Y_{C},Y_{B})Y_{A}=\mathbf{R}_{ABC}^{F}Y_{F},\qquad%
\mathbf{5.}\text{\hspace{3mm}}\nabla _{Y_{C}}\omega ^{B}=-\Gamma
_{FC}^{B}\omega ^{F},\medskip \\ 
\mathbf{6.} & [\mathbf{R}(Y_{C},Y_{B})X]\cdot \omega ^{B}\cdot \omega
^{C}=\left\{ \nabla _{Y_{C}}\nabla _{Y_{B}}X-\nabla _{Y_{B}}\nabla
_{Y_{C}}X-\right. \medskip \\ 
& \hspace{41mm}-\nabla _{\lbrack Y_{C},Y_{B}]}X\}\cdot \omega ^{B}\cdot
\omega ^{C},%
\end{array}%
\end{equation*}%
where \textquotedblright $\cdot $\textquotedblright\ represents the
tensorial product \textquotedblright $\otimes $\textquotedblright , we
deduce by a direct calculation that%
\begin{equation}
X_{:B:C}^{A}-X_{:C:B}^{A}=X^{F}\mathbf{R}_{FBC}^{A}-X_{:F}^{A}\mathbf{T}%
_{BC}^{F},  \label{ric}
\end{equation}%
where \textquotedblright $_{:D}$\textquotedblright\ represents one from the
local covariant derivatives \textquotedblright $_{/1}$\textquotedblright ,
\textquotedblright $_{|j}$\textquotedblright\ or \textquotedblright $%
|_{(j)}^{(1)}$\textquotedblright\ of the $\Gamma $-linear connection $\nabla 
$.

Taking into account that the indices $A,B,C,\ldots $ belong to the set $%
\left\{ 1,i,{\QATOP{(1)}{(i)}}\right\} $, by complicated local computations,
the identities (\ref{ric}) imply the required Ricci identities.
\end{proof}

Now, let us consider the \textit{canonical Liouville d-tensor field}%
\begin{equation*}
\mathbf{C}{=\mathbf{C}_{(1)}^{(i)}{\frac{\partial }{\partial y_{1}^{i}}=}}%
y_{1}^{i}{{\frac{\partial }{\partial y_{1}^{i}}.}}
\end{equation*}

\begin{definition}
The distinguished tensors defined by the local components 
\begin{equation}
\bar{D}_{(1)1}^{(i)}=\mathbf{C}_{(1)/1}^{(i)},\quad D_{(1)j}^{(i)}=\mathbf{C}%
_{(1)|j}^{(i)},\quad d_{(1)(j)}^{(i)(1)}=\mathbf{C}_{(1)}^{(i)}|_{(j)}^{(1)}
\label{defldef}
\end{equation}%
are called the \textbf{deflection d-tensors attached to the }$\Gamma $%
\textbf{-linear connection }$\nabla $\textbf{\ on the 1-jet space }$E=J^{1}(%
\mathbb{R},M)$.
\end{definition}

Taking into account the expressions of the local covariant derivatives of
the $\Gamma $-linear connection $\nabla $ given by (\ref{Nabla-Gamma}), by a
direct calculation, we find

\begin{proposition}
The deflection d-tensors of the $\Gamma $-linear connection $\nabla $ have
the expressions:%
\begin{equation}
\begin{array}{c}
\bar{D}_{(1)1}^{(i)}=-M_{(1)1}^{(i)}+G_{(1)(r)1}^{(i)(1)}y_{1}^{r},\quad
D_{(1)j}^{(i)}=-N_{(1)j}^{(i)}+L_{(1)(r)j}^{(i)(1)}y_{1}^{r},\bigskip  \\ 
\medskip d_{(1)(j)}^{(i)(1)}=\delta
_{j}^{i}+C_{(1)(r)(j)}^{(i)(1)(1)}y_{1}^{r}.%
\end{array}
\label{deflloc}
\end{equation}
\end{proposition}

In the sequel, applying the set $(v)$ of the Ricci identities to the
components of the canonical Liouville d-tensor field $\mathbf{C}$, we get

\begin{theorem}
The deflection d-tensors attached to the $\Gamma $-linear connection $\nabla 
$ on the 1-jet space $E=J^{1}(\mathbb{R},M)$ verify the following identities:%
\begin{equation*}
\left\{ 
\begin{array}{l}
\medskip \bar{D}%
_{(1)1|k}^{(i)}-D_{(1)k/1}^{(i)}=y_{1}^{r}R_{(1)(r)1k}^{(i)(1)}-\bar{D}%
_{(1)1}^{(i)}\bar{T}%
_{1k}^{1}-D_{(1)r}^{(i)}T_{1k}^{r}-d_{(1)(r)}^{(i)(1)}R_{(1)1k}^{(r)} \\ 
\medskip
D_{(1)j|k}^{(i)}-D_{(1)k|j}^{(i)}=y_{1}^{r}R_{(1)(r)jk}^{(i)(1)}-D_{(1)r}^{(i)}T_{jk}^{r}-d_{(1)(r)}^{(i)(1)}R_{(1)jk}^{(r)}
\\ 
\medskip \bar{D}%
_{(1)1}^{(i)}|_{(k)}^{(1)}-d_{(1)(k)/1}^{(i)(1)}=y_{1}^{r}P_{(1)(r)1(k)}^{(i)(1)\;\;(1)}-%
\bar{D}_{(1)1}^{(i)}\bar{C}%
_{1(k)}^{1(1)}-d_{(1)(r)}^{(i)(1)}P_{(1)1(k)}^{(r)\;\;(1)} \\ 
\medskip
D_{(1)j}^{(i)}|_{(k)}^{(1)}-d_{(1)(k)|j}^{(i)(1)}=y_{1}^{r}P_{(1)(r)j(k)}^{(i)(1)\;\;(1)}-D_{(1)r}^{(i)}C_{j(k)}^{r(1)}-d_{(1)(r)}^{(i)(1)}P_{(1)j(k)}^{(r)\;\;(1)}
\\ 
\medskip
d_{(1)(j)}^{(i)(1)}|_{(k)}^{(1)}-d_{(1)(k)}^{(i)(1)}|_{(j)}^{(1)}=y_{1}^{r}S_{(1)(r)(j)(k)}^{(i)(1)(1)(1)}-d_{(1)(r)}^{(i)(1)}S_{(1)(j)(k)}^{(r)(1)(1)}.%
\end{array}%
\right.
\end{equation*}
\end{theorem}

\textbf{Acknowledgements. }The present research was supported by Contract
with Sinoptix No. 8441/2009.

\textbf{Author's addresses:}\medskip

Mircea N{\scriptsize EAGU }and Emil S{\scriptsize TOICA}

University Transilvania of Bra\c{s}ov, Faculty of Mathematics and Informatics

Department of Algebra, Geometry and Differential Equations

B-dul Eroilor 29, BV 500036, Bra\c{s}ov, Romania.

\textit{E-mails}: mircea.neagu@unitbv.ro, e.stoica@unitbv.ro

\textit{Websites}: http://www.2collab.com/user:mirceaneagu

\hspace{16.5mm}http://cs.unitbv.ro/\symbol{126}geome

\end{document}